\documentclass[12pt,a4paper]{amsart}
\usepackage[all,dvips]{xy}
\usepackage{amssymb}
\usepackage{graphicx}

\newcommand{\Endo}{\rm{End}}
\newcommand{\Q}{\ensuremath{\mathbb{Q}}}
\newcommand{\Z}{\ensuremath{\mathbb{Z}}}

\newcommand{\C}{\ensuremath{\mathbb{C}}}
\newcommand{\Zp}{\ensuremath{\Z_{p}}}
\newcommand{\Zpu}{\ensuremath{\Zp^{\times}}}
\newcommand{\Zpl}{\ensuremath{\Z_{(p)}}}
\newcommand{\Zplu}{\ensuremath{\Zpl^{\times}}}

\newcommand{\Ztl}{\ensuremath{\Z_{(2)}}}

\newcommand{\floor}[1]{\left\lfloor#1\right\rfloor}
\newcommand{\gp}[2]{\genfrac{[}{]}{0pt}{}{#1}{#2}}
\newcommand{\gpv}[3]{\gp{#1}{#2}_{#3}}
\newcommand{\KUpz}{\ensuremath{KU_{(p)}^{0}(KU_{(p)})}}
\newcommand{\kupz}{\ensuremath{ku_{(p)}^{0}(ku_{(p)})}}
\newcommand{\kut}{\ensuremath{ku_{(2)}}}
\newcommand{\kutz}{\ensuremath{ku_{(2)}^{0}(ku_{(2)})}}
\newcommand{\kuz}{\ensuremath{ku^{0}(ku)}}
\newcommand{\KUp}{\ensuremath{KU_{(p)}}}
\newcommand{\kup}{\ensuremath{ku_{(p)}}}
\newcommand{\MUp}{\ensuremath{MU_{(p)}}}
\newcommand{\MUt}{\ensuremath{MU_{(2)}}}
\newcommand{\MUpz}{\ensuremath{MU_{(p)}^0(MU_{(p)})}}
\newcommand{\MUtz}{\ensuremath{MU_{(2)}^0(MU_{(2)})}}
\newcommand{\MUps}{\ensuremath{MU_{(p)}^*(MU_{(p)})}}
\newcommand{\hC}{\widehat{c}}
\newcommand{\qh}{\hat{q}}
\newcommand{\hphi}{\hat{\varphi}}
\newcommand{\thetab}{\bar{\theta}}
\newcommand{\pib}{\bar{\pi}}
\newcommand{\ab}{\bar{\alpha}}
\newcommand{\rb}{\bar{r}}
\newcommand{\Sb}{\bar{S}}

\theoremstyle{plain}
\newtheorem{theorem}[equation]{Theorem}
\newtheorem{proposition}[equation]{Proposition}

\newtheorem{lemma}[equation]{Lemma}

\theoremstyle{definition}
\newtheorem{definition}[equation]{Definition}

\begin{document}

\title{Infinite Sums of Adams Operations and Cobordism}
\author{Imma G\'alvez \and Sarah Whitehouse}                             
\address{Computing, Communications Technology and Mathematics, 
London Metropolitan University, Holloway Road, London N7 8DB, UK.}
\email{i.galvezicarrillo@londonmet.ac.uk} 
\address{Pure Mathematics, University of Sheffield, Sheffield S3 7RH, UK.} 
\email{s.whitehouse@sheffield.ac.uk}

\date{$11^{\text{th}}$ June 2004}

\begin{abstract}
In recent work, various algebras of stable degree zero operations in $p$-local
$K$-theory were described explicitly~\cite{ccw2}. The elements are certain 
infinite sums of Adams operations. Here we show how to make sense of 
the same expressions for $MU_{(p)}$ and
$BP$, thus identifying the ``Adams subalgebra'' of the algebras of operations.
We prove that the Adams subalgebra is the centre of the ring of degree
zero operations.
\end{abstract}


\maketitle

\section{Introduction}
\label{intro}

This paper builds on recent work on stable operations in $p$-local complex
$K$-theory~\cite{ccw2}.
That work described all operations in terms of Adams operations. In a
sense, this was done by specifying which infinite sums of Adams operations 
are allowed. Here we explore these expressions for $\MUp$ and $BP$.

We begin by briefly discussing the definition and basic properties of Adams operations
for these theories. These facts are well known, but we hope it
is helpful to gather them together here. The main results of~\cite{ccw2}
are then recalled and the analogues for $MU_{(p)}$ and $BP$ are proved, thus
identifying the ``Adams subalgebra'' of the algebras of operations for these
theories. 

The final sections are devoted to showing that this subalgebra 
is precisely the centre of the ring of stable degree zero operations.
This can be compared to old results of Novikov (for $MU$) and
Araki (for $BP$),~\cite{novikov67,araki}. 
These state that the Adams operations form the centre of
the \emph{group} of stable degree zero 
multiplicative operations.

\section{Adams Operations}
\label{sec2}
We begin by summarizing some standard results on stable cohomology
operations. 

\begin{proposition}\cite[Section 14]{boardman}
Let $E=MU$, $BP$ or $KU$. A stable multiplicative cohomology
operation $\theta:E\to E$ is uniquely determined by its value
$\theta(x_{E})$ on the orientation class $x_E\in
E^2(\mathbb{C}P^\infty)$.
\qed\end{proposition}

We remark that the same is true for the Adams summand of $p$-local 
$K$-theory, which we denote by $G$.

In order to have stable Adams operations, we will work $p$-locally, for
a prime $p$.

\begin{definition}
Let $E=\MUp$, $BP$, $\KUp$ or $G$.  A stable Adams operation 
$\Psi_E^{\alpha}$ is defined
for each ${\alpha}\in\Zplu$ as the unique 
multiplicative operation given on the
orientation class $x_E$ by
    $$
    \Psi_E^{\alpha}(x_E)=\frac{[{\alpha}]_E(x_E)}{{\alpha}}
    \in E^*(\mathbb{C}P^\infty)=E_*[[x_E]]\, ,
    $$
where $[{\alpha}]_E$ denotes the formal sum.
\end{definition}

\begin{lemma}\label{actcoeffs}
$\Psi_E^{\alpha}$ acts as multiplication by ${\alpha}^n$ on $E_{2n}$.
\end{lemma}

\begin{proof}
Let the $\log$ series of the formal group law corresponding to $E$
be given by $\log_E(x_E)=\sum_{i=0}^\infty m_ix_E^{i+1}$, where $m_i\in
E_{2i}\otimes\Q$. The operation $\Psi_E^{\alpha}$ transforms this 
into $\log_E\left((\Psi_E^{\alpha})^{-1}(x_E) \right)$, 
see~\cite[Appendix B]{ravorange}. This
is
    \begin{eqnarray*}
    \log_E\left( (\Psi_E^{\alpha})^{-1}(x_E) \right)
    &=&\log_E\big( \exp_E\frac{1}{{\alpha}}\log_E {\alpha}x_E\big)\\
    &=&\frac{\log_E{\alpha}x_E}{\alpha}
    =\sum_{i=0}^\infty {\alpha}^i m_i x_E^{i+1}.
    \end{eqnarray*}
Thus on coefficients, $m_i\mapsto {\alpha}^i m_i$. Since $E_*$ is torsion
free and since $\Psi_E^{\alpha}$ is multiplicative, 
this is sufficient to determine $\Psi_E^{\alpha}$ on the coefficient groups.
\hspace*{\fill}
\end{proof}

For $E=MU$, $BP$, $KU$ or $G$, the representation of the 
ring of operations given by their actions on coefficient groups 
is faithful:
    $$
    E^*(E)\hookrightarrow \Endo(E_*).
    $$
For $MU$ this result is due to Novikov~\cite[Lemma 5.2]{novikov67}; 
for $KU$, see~\cite[Theorem 2]{Johnson1984}.
The $BP$ and $G$ versions follow.
In particular, we see that Adams operations multiply as 
expected: $\Psi_E^\alpha\Psi_E^\beta=\Psi_E^{\alpha\beta}$. We also note that,
since $\Psi^\alpha_E$ is multiplicative, it is a grouplike element 
with respect to the coproduct of $E^*(E)$, 
$\Delta(\Psi^\alpha_E)= \Psi^\alpha_E\otimes\Psi^\alpha_E$.

If $e$ denotes the connective cover of $E$, the map $e\to E$ leads to a  
map $E^*(E)\to e^*(e)$. Thus we may define an Adams operation 
$\Psi_e^\alpha$ for $e$, with all the expected
properties, as the image of $\Psi_E^\alpha$ under this map.

\section{Operations in $p$-local $K$-theory}
\label{sec3}

Let $\KUp$ be the $p$-local periodic complex $K$-theory spectrum
and let $\kup$ be the connective version.
As above, for these theories we have stable Adams operations,
$\Psi_{KU}^\alpha$ and $\Psi_{ku}^\alpha$ respectively, for
each $\alpha$ in $\Zplu$. 

Let $p$ be an odd prime and let $q$ be primitive modulo $p^2$. This ensures
that $q$ is primitive modulo $p^r$ for all $r\geq 1$, which is to say that
the powers of $q$ are dense in the $p$-adic units $\Zpu$.
Clearly, the polynomial ring generated by the Adams operation
$\Psi_{KU}^q$ is a subring of the
ring of all stable degree zero $\KUp$ operations. 
This does not give all operations; it is countable and the operation algebra is 
known to be uncountable~\cite{ac}.

\begin{definition}
Let $p$ be an odd prime and let $q$ be primitive modulo $p^2$.
Let $E=\MUp$ or $\KUp$ and let $e=\MUp$ or $\kup$. 
We define two families of polynomials.
\begin{enumerate}
\item[1)] Let $\varphi^e_n= \prod_{i=0}^{n-1} (\Psi_e^q-q^i)\in 
\Zpl[\Psi_e^q]$. 
\item[2)] Let $\Phi^E_n=  \prod_{i=0}^{n-1} (\Psi_E^q-q_i) \in 
\Zpl[\Psi_E^q],$ where $q_i$
is the $i+1$-st term of the sequence
$$
  1,q, q^{-1}, q^2, q^{-2}, q^3, q^{-3}, q^4,\  \dots\ .
$$
\end{enumerate}
\end{definition}

For $p$ odd, the
following result of~\cite{ccw2} shows how to view the ring of all degree zero
stable $KU_{(p)}$ operations as a completion
of the polynomial
ring $\Zpl[\Psi_{KU}^q]$. It also gives a description for the
connective case, corresponding to a different completion of a
polynomial ring.
These results may be viewed as specifying which 
infinite sums of Adams operations are admissible.
\smallskip

\begin{theorem}~\cite[Theorems 2.2 and 6.2]{ccw2}\label{phi-basis}
\begin{enumerate}
\item[1)] Elements of $\kupz$ can be expressed uniquely as
    infinite sums
    $$\sum_{n \ge 0}  a_{n} \varphi^{ku}_{n} ,$$
    where $a_{n} \in \Z_{(p)}$.

\item[2)] Elements of $\KUpz$ can be expressed uniquely as
infinite sums
    $$\sum_{n \geq 0}  a_{n} \Phi^{KU}_{n} ,$$
    where $a_{n} \in \Z_{(p)}$.\qed
\end{enumerate}
\end{theorem}
Explicit product and coproduct formulas and many related
results can be found in~\cite{ccw2}. Analogues of these results 
hold for $p=2$, but they are more complicated. We give only the connective version.

\begin{definition}\label{zetadef}
    Define elements $\zeta^{ku}_{n}\in \kutz$, for $n\ge 0$, by
    \begin{eqnarray*}
        \zeta^{ku}_{2m+1}
        &=&
        (\Psi_{ku}^{-1}-1)\prod_{i=0}^{m-1}(\Psi_{ku}^{3}-3^{2i+1}),\\
        \zeta^{ku}_{2m}
        &=&
        \prod_{i=0}^{m-1}(\Psi_{ku}^{3}-3^{2i})
        +
        \sum_{i=1}^{m}
        \frac{\theta_{i}(3)\theta_{i}(3^{2m})}{2\theta_{i}(3^{2i})}
    \zeta^{ku}_{2m-2i+1},
    \end{eqnarray*}
where $\theta_r(X)=\prod_{i=0}^{r-1}(X-3^{2i})$.
\end{definition}

\begin{theorem}~\cite[Theorem 8.2]{ccw2}\label{zetabasis}
    The elements of $\kutz$ can be expressed uniquely as
    infinite sums
    $$
    \sum_{n\ge 0}a_{n}\zeta^{ku}_{n},
    $$
    where $a_{n}\in\Ztl$.\qed
\end{theorem}

\section{Infinite Sums of Adams Operations for $\MUp$}
\label{sec4}

In this section, we show that the infinite sums of Adams operations of
the previous section are also defined for $\MUp$. 

\begin{proposition}\label{infsumsMU}
Let $p$ be an odd prime.
The infinite sums $\sum_{n=0}^\infty a_n\varphi^{MU}_n\!$, where
$a_n\in\Zpl$, are well-defined operations in
$\MUp$-cohomology.
\end{proposition}

\begin{proof}
Clearly, finite sums of the $\varphi^{MU}_n$ are well-defined operations.
To see that the same is true for the infinite sums, it suffices
to show that $\varphi^{MU}_n\to 0$ as $n\to \infty$ in the usual
filtration topology of $MU_{(p)}^*(MU_{(p)})$.

This may be checked by considering the action on coefficient groups.
By Lemma~\ref{actcoeffs}, $\Psi_{MU}^q$ acts
as multiplication by $q^i$ on
$\pi_{2i}(MU_{(p)})$. It follows that $\varphi^{MU}_n$ acts as zero 
on $\pi_{2i}(MU_{(p)})$ for $i<n$.
\end{proof}

\begin{proposition}\label{kinject}
For any odd prime $p$,
the map 
    $$
    \kupz\to \MUpz
    $$ 
given by
    $$
    \sum_{n=0}^\infty a_n\varphi^{ku}_n
    \mapsto\sum_{n=0}^\infty a_n\varphi^{MU}_n
    $$ 
is an injective algebra map.
\end{proposition}

\begin{proof}
Consider $\sum_{n=0}^\infty
a_n\varphi^{MU}_n$ in $\MUpz$ and suppose that
$a_m\neq 0$, with $m$ chosen minimal. This
operation acts on $\pi_{2m}(\MUp)\neq 0$ as the
finite sum $\sum_{n=0}^m a_n\varphi^{MU}_n=a_m\varphi^{MU}_m$ and
 thus as multiplication by $a_m\prod_{i=0}^{m-1}(q^m-q^i)\neq
0$. Thus the operation is non-trivial. So the map is injective.

It is easy to see that we have an algebra map: 
the product of two 
infinite sums is determined in both the source and the target by the products
of Adams operations. (See~\cite[Prop.\,2.7]{ccw2} for an 
explicit formula.)
\end{proof}

The source $\kupz$ is a completed bialgebra.
We note that the injective map above also respects the
coproduct, since the coproduct of a general infinite sum is determined 
by the fact that the Adams operations are
group-like; see~\cite[Prop.\,2.9]{ccw2} for an explicit formula.
So the image of $\kupz$ is a subbialgebra of 
$\MUps$.

There is also an injective algebra map $\KUpz\to \MUpz$. 
This is because the injection of Proposition~\ref{kinject} 
can be composed with
the inclusion of  bialgebras
    $$
    \KUpz\hookrightarrow \kupz,
    $$ 
resulting from the covering map $\kup\to \KUp$. 
An explicit formula for this map is
given in~\cite[Prop.\,7.1]{ccw2}, expressing each $\Phi^{KU}_n$ 
in terms of the $\varphi^{ku}_m$. 

Again the injection is well-behaved with respect to the coproduct.
This time we have a conjugation map $\chi$ and this is also respected
by the inclusion, as it is determined for both 
$E=\KUp$ and $E=\MUp$ by 
$\chi\Psi_E^q=\Psi_E^{q^{-1}}$.
\smallskip

So far the discussion in this section has all been for odd primes.
However, it is easy to see that the analogues for $p=2$ also hold. 
One defines
$\zeta_n^{MU}\in \Ztl[\Psi_{MU}^3, \Psi_{MU}^{-1}]$ 
in the obvious way, just as in Definition~\ref{zetadef}.

\begin{proposition}\label{ktinject}
\begin{enumerate}
\item The infinite sums $\sum_{n=0}^\infty a_n\zeta^{MU}_n$, where
$a_n\in\Ztl$, are well-defined operations in
$\MUt$-cohomology.
\item The map $\kutz\to \MUtz$ given by
    $$
    \sum_{n=0}^\infty a_n\zeta^{ku}_n
    \mapsto\sum_{n=0}^\infty a_n\zeta^{MU}_n
    $$ 
is an injective algebra map.
\end{enumerate}
\end{proposition}

\begin{proof}
The key properties of $\zeta^{ku}_n$ are that it acts as zero on coefficient
groups $\pi_{2i}(\kut)$, for all $i< n$, and its action is non-zero
on $\pi_{2n}(\kut)$. (See~\cite[Section 8]{ccw2} 
for a proof and related formulas.)
Using this, one may check that the infinite sums are
well-defined $\MUt$ operations just as in Proposition~\ref{infsumsMU} and that
we have the claimed injection of algebras just as in Proposition~\ref{kinject}.
\end{proof}

\section{Infinite Sums of Adams Operations for $BP$}
\label{sec5}

We record here the $BP$ analogues of the results of the preceding section, 
omitting proofs since these are easy modifications of the $\MUp$
versions.

For $p$ an odd prime, the spectra $\KUp$ and~$\kup$ each
split into $p-1$ copies of spectra which we denote by $G$
and~$g$, respectively%
\footnote{The notations $E(1)$ and $e(1)$ and $L$ and $l$ are also used.}. 
As before we choose $q$ primitive modulo $p^2$ and now we let 
$\qh=q^{p-1}$. Thus the powers of $\qh$ are dense in $1+p\Zp$.
Now we recall from~\cite{ccw2} the description of degree zero 
stable operations for
$g$.

\begin{definition}
Define $\hphi^g_n \in \Zpl[\Psi_g^q]\subset g^0(g)$, for $n\geq 0$, by
     $$
     \hphi^g_n=\prod_{i=0}^{n-1}(\Psi_g^{q}-\qh^{i}) .
     $$
Also define $\hphi^{BP}_n \in \Zpl[\Psi_{BP}^q]\subset BP^0(BP)$
in the same way.     
\end{definition}
\smallskip

\begin{theorem}~\cite[Theorem 2.2]{Lellmann1984},~\cite[Theorem 4.4]{ccw2}\label{g-ops}
Elements of $g^{0}(g)$ can be expressed uniquely as infinite sums
    $$\sum_{n \geq 0}  a_{n} \hphi^g_{n} ,$$
    where $a_{n} \in \Z_{(p)}$.\qed
\end{theorem}

\begin{proposition}\label{ginject}
Let $p$ be an odd prime.
\begin{enumerate}
\item The infinite sums $\sum_{n=0}^\infty a_n\hphi^{BP}_n$, where
$a_n\in\Zpl$, are well-defined operations in
$BP$-cohomology.
\item The map $\iota : g^0(g)\to BP^0(BP)$ given by
    $$
    \sum_{n=0}^\infty a_n\hphi^g_n
    \mapsto\sum_{n=0}^\infty a_n\hphi^{BP}_n
    $$ 
is an injective algebra map.\hfill\qed
\end{enumerate}
\end{proposition}

The same remarks as in the preceding section about the coproduct and 
the comparison
with the periodic operations may be made here.

Now let $p=2$. 
Again we define
$\zeta_n^{BP}\in \Ztl[\Psi_{BP}^3, \Psi_{BP}^{-1}]$ 
just as in Definition~\ref{zetadef}.

\begin{proposition}\label{kubpinject}
Let $p=2$.
\begin{enumerate}
\item The infinite sums $\sum_{n=0}^\infty a_n\zeta^{BP}_n$, where
$a_n\in\Ztl$, are well-defined operations in
$BP$-cohomology.
\item The map $\iota : \kutz\to BP^0(BP)$ given by
    $$
    \sum_{n=0}^\infty a_n\zeta^{ku}_n
    \mapsto\sum_{n=0}^\infty a_n\zeta^{BP}_n
    $$ 
is an injective algebra map.\hfill\qed
\end{enumerate}
\end{proposition}

\section{Diagonal Operations}
\label{sec6}

We will need to discuss operations which
act diagonally on coefficient groups and so we introduce the 
following notation.

\begin{definition} Let $E=MU$, $\MUp$ or $BP$. We write $D_{E}$ 
for the subring of
$E^0(E)$ consisting of operations whose action on each coefficient group
$E_{2i}$ is multiplication by an element $\lambda_i$ of the ground ring, (that is
an element of $\Z$ for $E=MU$, an element of $\Zpl$ for $E=\MUp$ or $E=BP$). 
\end{definition}

We next recall some results from~\cite{novikov67} on $MU$ operations. 
For each finite non-decreasing 
sequence of positive integers, $\alpha=(\alpha_1, \alpha_2, \alpha_3, \dots)$,
there is an operation $\Sb_\alpha\in MU^*(MU)$, 
of degree $2|\alpha|=2\sum_i \alpha_i$.
Here $\Sb_\alpha=\phi(\bar{\sigma}_\alpha)$, where $\phi:MU^*(BU)\to MU^*(MU)$
is the Thom isomorphism, $\sigma_\alpha\in MU^*(BU)$ is the Conner-Floyd characteristic class
associated to $\alpha$
and $\bar{\sigma}_\alpha$ is defined by $\bar{\sigma}_\alpha(\xi)=\sigma_\alpha(-\xi)$.

\begin{proposition}\label{diag}
Let $E=MU$, $\MUp$ or $BP$. Then $Z\left(E^0(E)\right)=D_E$.
\end{proposition}

\begin{proof}
Clearly, $D_E$ is contained in the centre $Z\left(E^0(E)\right)$. 
To prove the reverse inclusion, we first let $E=MU$. Since $E_*$ is 
torsion free, an operation $\theta$ is zero
if and only if $\theta_*\otimes 1_{\Q}$ is zero on 
$E_*\otimes\Q=\Q[\C P^1, \C P^2, \dots ]$. So we can work
with monomials in the classes $[\C P^n]$. We first consider the 
operation
$\Sb_{(m)}\in E^{2m}(E)$. We let $\mathcal{S}_{(m)}=[\C P^m]\Sb_{(m)}\in E^0(E)$. This
operation
has the property that it acts as zero on any decomposable class of degree $2m$
and as multiplication by $(m+1)$ on $[\C P^m]$. 
Let $\theta\in Z(E^0(E))$.
The relation $\mathcal{S}_{(m)}\theta=\theta \mathcal{S}_{(m)}$ then tells us that
$\theta_*[\C P^m]=\lambda_m[\C P^m]$ for some $\lambda_m\in \Z$.

More generally, for $\alpha=(\alpha_1,\alpha_2,\dots,\alpha_r)$, 
we have a Landweber-Novikov operation
$\Sb_\alpha\in E^{2|\alpha|}(E)$. By~\cite[Lemma 5.5]{novikov67},  
$(\Sb_\alpha)_*([\C P^n])=\lambda_\alpha[\C P^{n-|\alpha|}]$, 
where $\lambda_\alpha\in\Z$ and $\lambda_\alpha\neq 0$.
Let $\C P^\alpha=\C P^{\alpha_1}\times\C P^{\alpha_2}\times\dots\times\C P^{\alpha_r}$ 
and let $\mathcal{S}_\alpha=[\C P^\alpha]\Sb_\alpha\in E^0(E)$.
Then the relation 
$(\mathcal{S}_{\alpha}\theta)_*[\C P^{|\alpha|}]
=(\theta\mathcal{S}_{\alpha})_*[\C P^{|\alpha|}]$ tells us that $\theta$ acts
on $[\C P^\alpha]$ as multiplication by $\lambda_{|\alpha|} \in \Z$.
Thus $\theta\in D_{E}$.

The argument is just as above for $\MUp$ and that for $BP$ is similar, 
but here we need only work with
monomials in the $[\C P^{p^i-1}]$, since 
$BP_*\otimes\Q=\Q[\C P^{p^i-1}\,|\, i\geq 1]$. 
\end{proof}

We remark that, if one considers the algebra of all operations, 
rather than only degree zero ones, then the centre consists of just 
the constant operations.

\section{Congruences}
\label{sec7}

In this section $p$ is an odd prime and we compare two families of
congruences, one related to the connective Adams summand $g$ and the
other to $BP$. The results will be used in the next
section to give a new description of the centre of $BP^0(BP)$.

Clearly, we may identify the ring of diagonal operations
$D_{BP}$ with a subring of the infinite
direct product $\prod_{k\geq 0}\Zpl$. In what follows we will often use 
this identification without further comment.

We begin by noting the analogue for the connective Adams summand $g$
of~\cite[Theorem 11]{ccw}. This characterizes $g^0(g)$ as a subring 
of the infinite direct product $\prod_{k\geq 0}\Zpl$ by a family 
of congruences which involve Gaussian polynomials. 

\begin{definition}
The Gaussian polynomial $\gp{n}{i}\in \Z[t]$ is defined,
for non-negative integers $n$ and $i$, by
   $$
   \gp{n}{i}=\prod_{k=0}^{i-1} \frac{1-t^{n-k}}{1-t^{i-k}}.
   $$
Also let $\gpv{n}{i}{a}$ denote the
value of this polynomial at $t=a$.
\end{definition}

\begin{theorem}\label{gcong}
If $\phi\in g^0(g)$ acts on $g_{2(p-1)i}=\Zpl$ as 
multiplication by $\mu_i\in\Zpl$, then 
    $$
    \sum_{i=0}^n (-1)^{n-i} \qh^{\binom{n-i}{2}}\gpv{n}{i}{\qh} \mu_i 
    \equiv 0\quad\bmod{p^{\delta_{p}(n)}}, 
    $$
for all $n\geq 0$, where $\delta_p(n)=n+\nu_p(n!)$. 
Moreover every sequence satisfying these congruences arises from a 
unique stable
operation.
\end{theorem}

\begin{proof}
A $\Zpl$-basis for $G_0(g)$ is given in~\cite[Prop.\,4.2]{ccw2}.
Using~\cite[Prop.\,8]{ccw} this may be expressed in terms
of Gaussian polynomials. The result
follows from the fact that $g^0(g)$ is the $\Zpl$-linear dual of $G_0(g)$.
\end{proof}

\begin{definition}\label{coeffs}
Let
    $$
    C_{n,i}= \frac{ (-1)^{n-i} \qh^{\binom{n-i}{2}}\gpv{n}{i}{\qh} }
                    { p^{\delta_{p}(n)}  } \  \in\Q,
    $$
and let $C_n=(C_{n,i})_{i\geq 0}$, a sequence of rational numbers.
\end{definition}

We note that $C_{n,n}=p^{-\delta_p(n)}$ and $C_{n,i}=0$ when $i>n$.
Also note  that $\delta_p(p^r)=1+p+p^2+\cdots +p^r$.
The $n$th congruence of Theorem~\ref{gcong} may be 
expressed as $\sum_{i=0}^n C_{n,i}\mu_i \in \Zpl$
or, adopting vector notation, $C_n\cdot \mu \in \Zpl$.

The rest of this section will be devoted to
comparing the congruences of Theorem~\ref{gcong} 
with a system of congruences characterizing $D_{BP}$.

Next we recall that, using duality, an operation $\theta:BP\to BP$ 
is determined by the left $BP_*$-module map $\thetab=\langle \theta,\ \rangle
:BP_*(BP)\to BP_*$, where $\langle\ ,\ \rangle$ is the Kronecker pairing.
The action of the operation on coefficient groups $\theta_*:BP_*\to BP_*$ 
is obtained by composition with the right unit map: $\theta_*=\thetab\eta_R$.
(See, for example,~\cite[11.22]{boardman}.)
\begin{figure}[h]
\begin{center}
\leavevmode
\includegraphics[240, 647][372, 714]{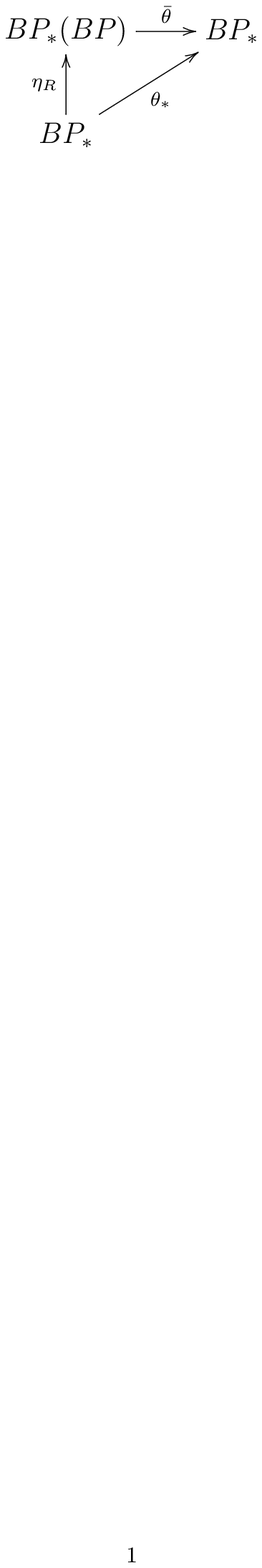}\\
\end{center}
\end{figure}

We will adopt the $BP$ notation of~\cite{ravorange}. Thus
$BP_*=\Zpl[v_1,v_2,\dots]$, where the $v_i$ are Araki's generators,
$|v_i|=2p^i-2$, and 
$BP_*(BP)=BP_*[t_1,t_2, \dots]$, where $|t_i|=2p^i-2$. We also use
the elements $l_i\in BP_*\otimes \Q$, related recursively to the
$v_i$ by the formula of Araki. As usual, because all objects are torsion
free, one may calculate rationally. We will do this implicitly, so that,
for example, we write $\thetab(l_i)$ to mean $(\thetab\otimes 1_{\Q})(l_i)$.

We suppose that we have $\theta\in D_{BP}$, acting on $BP_{2(p-1)i}$
as multiplication by $\mu_i \in \Zpl$.
(So, in the notation of Sect.~\ref{sec6}, $\mu=(\lambda_0, \lambda_{p-1},
\lambda_{2(p-1)}, \dots)$.)
By consideration of the relation $\theta_*=\thetab\eta_R$, 
one may calculate recursively $\thetab(x)$ for $x\in BP_*(BP)$ 
in terms of the sequence $\mu=(\mu_i)_{i\geq 0}$. Since, in $\thetab(x)$, 
the coefficient of each monomial in the generators $v_i$ must lie 
in $\Zpl$, one finds congruences imposed on the terms of this sequence.

We now give more details of this system of congruences. Let 
$\alpha, \beta, \gamma, \delta, \epsilon$ be finite sequences of 
non-negative integers. If $\alpha=(\alpha_1, \alpha_2, \dots, \alpha_m)$ 
then $v^\alpha$ denotes the monomial 
$v_1^{\alpha_1}v_2^{\alpha_2}\dots
v_m^{\alpha_m}$ and let $|\!|\alpha|\!|=\sum_{i=1}^m \alpha_i(1+p+\cdots +p^{i-1})
=\frac{|v^\alpha|}{2(p-1)}$.

Let $E_{\beta,\gamma}^\alpha
\in \Zpl$ be determined by 
$\eta_R(v^\alpha)=
\sum_{\beta,\gamma}E_{\beta,\gamma}^\alpha v^\beta t^\gamma$. 
Then, since $\thetab\eta_R(v^\alpha)=\theta_*(v^\alpha)$, we find
that
    \begin{equation}\label{explicitbpcongs}
    \sum_{\beta+\delta=\epsilon} \sum_{\gamma} E_{\beta, \gamma}^\alpha D_{\delta}^\gamma
    (\mu)=\begin{cases}
            \mu_{|\!|\alpha|\!|}&\mbox{if $\alpha=\epsilon$,}\cr
            0&\mbox{otherwise,}
            \end{cases}
    \end{equation}
where $\thetab(t^\gamma)=\sum_{\delta} D_\delta^\gamma(\mu)v^\delta$. Now
(\ref{explicitbpcongs}) determines $D_\delta^\gamma(\mu)$ recursively as a
finite rational linear combination of the $\mu_i$. On the other hand, we have $D_\delta^\gamma(\mu)\in \Zpl$,
for all $\gamma$ and $\delta$, and these integrality conditions 
form our family of congruences.

For example,  
$D_{(1)}^{(1)}(\mu)=a(\frac{\mu_1-\mu_0}{p})$, 
where $a\in \Zplu$.
So $\frac{\mu_1- \mu_0}{p}\in \Zpl$. Notice that this 
is the congruence $C_1\cdot\mu\in\Zpl$  of Theorem~\ref{gcong}.

In order to compare
the congruences of Theorem~\ref{gcong} and the $BP$ ones, we need
some further notation.

\begin{definition}
Let $\mu$  denote an element in $\prod_{k\geq 0}\Zpl$ and define the 
following subsets of $\prod_{k\geq 0}\Zpl$.
    \begin{eqnarray*}
    S_\infty^g&=&\left\{ \mu\,|\, C_i\cdot\mu\in\Zpl\
                        \mbox{for all $0\leq i$} \right\},\\
    S_n^g&=&\left\{ \mu\,|\, C_i\cdot\mu\in\Zpl\
                        \mbox{for $0\leq i\leq n$} \right\},\\
    S_\infty^{BP}&=&\left\{\mu\,|\, D_\delta^\gamma(\mu)\in\Zpl 
    \mbox{ for all $\gamma, \delta$} \right\},\\
    S_n^{BP}&=&\left\{\mu\,|\, \exists\ \mu'\in S_\infty^{BP}\ \mbox{with}\ 
    \mu_i=\mu'_i\ \mbox{for}\ 0\leq i\leq n \right\}.
    \end{eqnarray*}
\end{definition}

Thus $g^0(g)\cong S_\infty^g$ by Theorem~\ref{gcong} and $D_{BP}\cong S_\infty^{BP}$
by the preceding discussion. Also, clearly
    \begin{eqnarray*}
    S_0^g\supset S_1^g\supset\dots\supset S_n^g\supset S_{n+1}^g\supset\dots 
\supset S_\infty^g&=&\bigcap\limits_{n=0}^\infty S_n^g,\\
    S_0^{BP}\supset S_1^{BP}\supset\dots\supset S_n^{BP}\supset S_{n+1}^{BP}\supset\dots 
       \supset S_\infty^{BP}&=&\bigcap\limits_{n=0}^\infty S_n^{BP}.
    \end{eqnarray*}

\begin{lemma}\label{congincl}
For all $n\geq 0$, $S_n^g\subseteq S_n^{BP}$.
\end{lemma}

\begin{proof}
Let $\mu\in S_n^g$. Given the form of the congruences of Theorem~\ref{gcong}, 
it is easy to see that
there is some $\mu'\in S_\infty^g\cong g^0(g)$ such that $\mu_i=\mu'_i$ for 
$0\leq i\leq n$.
Then, using the injection $\iota$ of Proposition~\ref{ginject}, we can 
find an operation $\theta\in BP^0(BP)$
inducing $\mu'$. Thus $\mu'\in S_\infty^{BP}$. But then $\mu\in S_n^{BP}$.
\end{proof}

In due course, we will establish the reverse inclusion, using induction on $n$.
Note that the congruence $C_{n}\cdot \mu \in \Zpl$ has the form
    $$
    p^{-\delta_p(n)} (\mu_{n} + y ) \in \Zpl,
    $$
where $y$ is a $\Zpl$ linear combination of $\mu_{n-1}, \dots, \mu_{0}$.
We will see that it is enough to find a $BP$ congruence
like this for each $n$. Although the following lemma is rather elementary, for convenience 
we record the precise formulation we need, leaving the proof to the reader.

\begin{lemma}\label{reducetospecialcong}
For each $r\geq 0$, let $c_r=(c_{r,i})_{i\geq 0}$ be a sequence 
of rational numbers such that 
\begin{eqnarray*}
c_{r,i}&\in & p^{-\delta_p(r)}\Zpl, \quad\mbox{for $0\leq i\leq r-1$},\\ 
c_{r,r}&\in & p^{-\delta_p(r)}\Zplu,\\
c_{r,i}&= &0, \quad\mbox{for $i>r$}. 
\end{eqnarray*}
Let $\hC_r$ be another sequence satisfying the 
same conditions as $c_r$.

\noindent Let $n\geq 0$ and let
    \begin{eqnarray*}
    S_n&=&\left\{ \mu \in \prod_{k\geq 0}\Zpl\ \Big|
            \ c_r\cdot \mu \in \Zpl \mbox{ for } 0\leq r\leq n\right\},\\
    T_n&=&\left\{\mu \in \prod_{k\geq 0}\Zpl\ \Big|
            \ c_r\cdot \mu \in \Zpl \mbox{ for } 0\leq r\leq n-1,\ 
            \hC_n\cdot \mu \in \Zpl 
     \right\}.
    \end{eqnarray*}
If $S_n\subseteq T_n$ then $S_n=T_n$.\qed
\end{lemma}

Of course, the $C_r$ of Definition~\ref{coeffs} satisfy the
hypotheses of the $c_r$ in the lemma. 

\begin{definition}
Let $\theta\in D_{BP}$, with $\theta_*=\mu\in S_\infty^{BP}$.
Define $V_{\mu}:BP_{*}(BP)\to \Zpl$ as the composite
$V_{\mu}=\pi\thetab$, where $\pi:BP_*\to\Zpl$ is the algebra
map determined by $\pi(v_1)=1$ and $\pi(v_i)=0$ for all $i>1$. 
\end{definition}

Thus, for $x\in BP_*(BP)$, $V_{\mu}(x)$ is a finite rational linear 
combination of the
$\mu_i$ and the fact that it lies in $\Zpl$ is one of our $BP$
congruences. We will see that it is enough to consider only
those $BP$ congruences arising in this way.

\begin{lemma}\label{Vprod}
Let $x,y\in BP_*(BP)$. Suppose $V_{\mu}(x)=\sum_{i=0}^r a_i\mu_i$ and  $V_{\mu}(y)=\sum_{j=0}^s b_j\mu_j$. Then 
$V_{\mu}(xy)=\sum_{i=0}^r \sum_{j=0}^s a_ib_j\mu_{i+j}$.
\end{lemma}

\begin{proof}
This is easy to check, using that $\eta_R$ is a ring homomorphism, that
$\thetab$ is a left $BP_*$-module homomorphism 
and that $\thetab\eta_R=\theta_*=\mu$.
\end{proof}

We recall the standard notation $\pi_n=p-p^{p^n}$. We write,
for $n\geq 1$, 
$\pib_n=\frac{\pi_n}{p}=1-p^{p^n-1} \in \Zplu$ and 
$\ab_n=\prod_{i=1}^n \pib_i\in \Zplu$. 

\begin{lemma}\label{lv1}
$V_{\mu}(1)=\mu_0$ and for $i\geq 0$,  
    \begin{equation}\label{timage}
    V_{\mu}(t_{i+1})= p^{-(i+1)}\ab_{i+1}^{-1}   \mu_{\delta_p(p^i)}
    -\sum_{k=1}^{i+1} p^{-k}\ab_k^{-1} V_{\mu}(t_{i+1-k}^{p^k}).
    \end{equation}
\end{lemma}

\begin{proof}
It is easily checked using induction and the Araki formula
that, for $k\geq 1$, the coefficient of $v_1^{\delta_p(p^{k-1})}$ in $l_k$ 
is $p^{-k}\ab_k^{-1}\in \Q$.
We recall that
$\eta_R(l_{i+1})=\sum_{k=0}^{i+1} l_k t_{i+1-k}^{p^k}$ and thus
    $$
    \mu_{\delta_p(p^i)}l_{i+1}=\theta_*(l_{i+1})
    =\thetab\eta_R(l_{i+1})=\sum_{k=0}^{i+1} l_k \thetab(t_{i+1-k}^{p^k}).
    $$
Now~(\ref{timage}) follows by equating coefficients of 
$v_1^{\delta_p(p^i)}$. 
\end{proof}

Now we produce our special $BP$ congruence.

\begin{proposition}\label{specialcong}
For each $n\geq 0$, there is an element $d_n\in BP_*(BP)$ such that
$V_{\mu}(d_n)=\sum_{j=0}^n d_{n,j}\mu_j$,
where $d_{n,j}\in p^{-\delta_p(n)}\Zpl$ for $0\leq j\leq n-1$ and 
$d_{n,n}\in p^{-\delta_p(n)}\Zplu$.
\end{proposition}

\begin{proof}
Firstly, we note that it is enough to prove this for $n=p^i$ 
for each $i\geq 0$.
For suppose that we have found $d_{p^i}$, for $i\geq 0$, then, writing 
$N=\sum_{k=0}^M a_kp^k$,
where $0\leq a_k\leq p-1$, it is straightforward to check using 
Lemma~\ref{Vprod} 
that we may put $d_N=\prod_{k=0}^M d_{p^k}^{a_k}$.

Now we will construct $d_{p^i}$ inductively. We adopt the induction hypothesis
that such an element exists and moreover can be chosen of the form 
$t_{i+1}+pr_{i+1}$, for some $r_{i+1}\in BP_*(BP)$.
For $i=0$,
we may take $d_1=t_1$, since $V_\mu(t_1)=p^{-1}\ab_1^{-1}(\mu_1-\mu_0)$. 
Now we assume that we have constructed 
$d_{p^j}$ for $0\leq j<i$ and we explain how to construct $d_{p^i}$.

For the remainder of this proof, for $a$ and $b$ 
finite rational linear combinations of the $\mu_k$,
we write $a\sim b$ to mean that $a$ and $b$ differ
by some element of the form $a_0\mu_0+\cdots +a_{p^i}\mu_{p^i}$,
where $a_k\in p^{-\delta_p(p^i)}\Zpl$ for $0\leq k\leq p^i-1$ and
$a_{p^i}\in p^{-\delta_p(p^i)+1}\Zpl$. 
Note that, if $a\sim b$ and
if $a$ satisfies the conditions required for $V_{\mu}(d_{p^i})$, 
then so does $b$.

Using Lemma~\ref{Vprod}, it is not hard to see that 
we can correct for the first term on the right-hand side of 
the equation~(\ref{timage}).
In fact,
we can find
$\rb_0\in BP_*(BP)$, a $\Zpl$ linear combination of
$d_1^{\delta_p(p^i)}, \dots, d_1^{p^i+1}$, such that
    \begin{equation}\label{equiv1}
    V_{\mu}(t_{i+1}-p\rb_0)\sim 
    -\sum_{k=1}^{i} p^{-k}\ab_k^{-1} V_{\mu}(t_{i+1-k}^{p^k}).
    \end{equation}
    
Let $1\leq k\leq i$. By our induction hypothesis, 
$d_{p^{i-k}}=t_{i-k+1}+pr_{i-k+1}$, for some $r_{i-k+1}\in BP_*(BP)$. 
It follows that
    \begin{equation}\label{equiv2}
    d_{p^{i-k}}^{p^k}=(t_{i-k+1}+pr_{i-k+1})^{p^k}
    =t_{i-k+1}^{p^k} + p^{k+1}\rb_{i-k+1},
    \end{equation}
for some $\rb_{i-k+1}\in BP_*(BP)$.

Now we set 
    $$
    d_{p^i}=t_{i+1}-p\rb_0-p\sum_{k=1}^{i} \ab_k^{-1}\rb_{i-k+1}\ \in BP_*(BP).
    $$
Clearly $d_{p^i}$ has the form required by the induction hypothesis.
Using~(\ref{equiv1}) and~(\ref{equiv2}), a calculation shows that
    $$
      V_{\mu}(d_{p^i})
                \sim -\sum_{k=1}^{i} p^{-k}\ab_k^{-1} 
         V_{\mu}(d_{p^{i-k}}^{p^k}).
    $$

Now, using the induction
hypothesis and Lemma~\ref{Vprod}, one checks that 
$V_{\mu}(d_{p^i})$ has the
required form.
\end{proof}

Now fix a choice of $d_n\in BP_*(BP)$ as in Proposition~\ref{specialcong} and let $C_n^{BP}=(d_{n,i})_{i\geq 0}$.
Then the $BP$ congruence $V_{\mu}(d_n)\in\Zpl$ may be written 
$C_n^{BP}\cdot \mu \in \Zpl$.

\begin{proposition}\label{revcongincl}
For all $n\geq 0$, $S_n^{BP}= S_n^g$.
\end{proposition}

\begin{proof}
This will be proved by induction on $n$.
It is true for $n=0$, since $S_0^{BP}=S_0^g=\prod_{k\geq 0}\Zpl$. 
Now we assume that $S_{n-1}^{BP}= S_{n-1}^{g}$.
Let 
   $$
   T_n^{BP}=\left\{\mu\in\prod_{k\geq 0}\Zpl\, \Big|\, C_r\cdot\mu\in \Zpl\ 
      \mbox{for $0\leq r\leq n-1$, $C_n^{BP}\cdot\mu\in \Zpl$} \right\}.
   $$ 
By Lemma~\ref{congincl}, $S_n^g\subseteq S_n^{BP}$. Now if $\mu\in S_n^{BP}$,
then $\mu\in S_{n-1}^{BP}=S_{n-1}^g$, by the induction hypothesis, and so
$C_r\cdot\mu\in\Zpl$ for $0\leq r\leq n-1$. Also, if $\mu\in S_n^{BP}$, there
is $\mu'\in S_\infty^{BP}$ such that $\mu_i=\mu'_i$ for $0\leq i\leq n$.
By Proposition~\ref{specialcong}, $C_n^{BP}\cdot \mu'\in\Zpl$. But $C_n^{BP}\cdot \mu'=C_n^{BP}\cdot \mu$.
So $S_n^{BP} \subseteq T_n^{BP}$.
Thus we have
    $$
    S_n^g\subseteq S_n^{BP} \subseteq T_n^{BP}.
    $$
Applying 
Lemma~\ref{reducetospecialcong} with $c_r=C_r$ and $\hC_r=C_r^{BP}$
gives $S_n^g=T_n^{BP}$. So $S_n^{BP}= S_n^g$. 
\end{proof}

\section{The Centre of the Ring of Degree Zero Operations}
\label{sec8}

\begin{definition}
\begin{enumerate}
\item We define the Adams subalgebra of $\MUp$ operations, 
$A_{\MUp}$, by
    \begin{eqnarray*}
    A_{\MUp}&=&\left\{ \sum_{n=0}^\infty a_n\varphi^{MU}_n\,\Big|\, 
                     a_n\in\Zpl\right\}\subset \MUpz,
                     \quad\mbox{for $p$ odd,}\\
    A_{\MUt}&=&\left\{ \sum_{n=0}^\infty a_n\zeta^{MU}_n\,\Big|\, 
                     a_n\in\Ztl\right\}\subset \MUtz.
    \end{eqnarray*}
\item We define the Adams subalgebra of $BP$ operations, $A_{BP}$, by
    \begin{eqnarray*}
    A_{BP}&=&\left\{ \sum_{n=0}^\infty a_n\hphi^{BP}_n\,\Big|\, 
                     a_n\in\Zpl\right\}\subset BP^0(BP), \quad\mbox{for $p$ odd,}\\
    A_{BP}&=&\left\{ \sum_{n=0}^\infty a_n\zeta^{BP}_n\,\Big|\, 
                     a_n\in\Ztl\right\}\subset BP^0(BP), \quad\mbox{for $p=2$}.
    \end{eqnarray*}
\end{enumerate}
\end{definition}

Each Adams subalgebra is the image of an injective algebra map constructed
in Section 4 or 5 and it follows that there are isomorphisms:
    \begin{eqnarray*}
    &A_{\MUp}\cong \kupz, \qquad &A_{BP}\cong g^0(g), \qquad\mbox{for $p$ odd},\\
    &A_{\MUt}\cong \kutz, \qquad &A_{BP}\cong \kutz, \qquad\mbox{for $p=2$.}
    \end{eqnarray*}

In Section 6, we identified the centre of the ring of all degree
zero stable operations with the ring of diagonal operations.
In this section we show that this is precisely
the Adams subalgebra. Clearly, the Adams subalgebra is contained 
in the diagonal operations; we need to see that this inclusion
is in fact an equality.
We begin with $BP$ at odd primes, since this
case follows directly from the results of the previous section.

\begin{theorem}\label{BPoddcentre}
Let $p$ be an odd prime. The centre of the ring of 
degree zero $BP$ operations is the 
Adams subalgebra:
$$
Z\left(BP^0(BP)\right)=
A_{BP}.
$$
\end{theorem}

\begin{proof}
By Proposition~\ref{diag}, the centre is equal to the
ring of diagonal operations, so it is enough to show that
$D_{BP}=A_{BP}$.
By Proposition~\ref{ginject} and Theorem~\ref{gcong}, 
$A_{BP}\cong g^0(g)\cong S_\infty^g$. Also $D_{BP}\cong S_\infty^{BP}$.
Here the isomorphisms
are given by sending operations to their actions on coefficients.
By Proposition~\ref{revcongincl}, $S_n^{BP}= S_n^g$ 
for all $n$ and so 
$S_\infty^{BP}= S_\infty^g$. 
So we have a commutative diagram
\vspace{-0.5cm}
\begin{figure}[h]
\begin{center}
\leavevmode
\includegraphics[273, 653][310, 710]{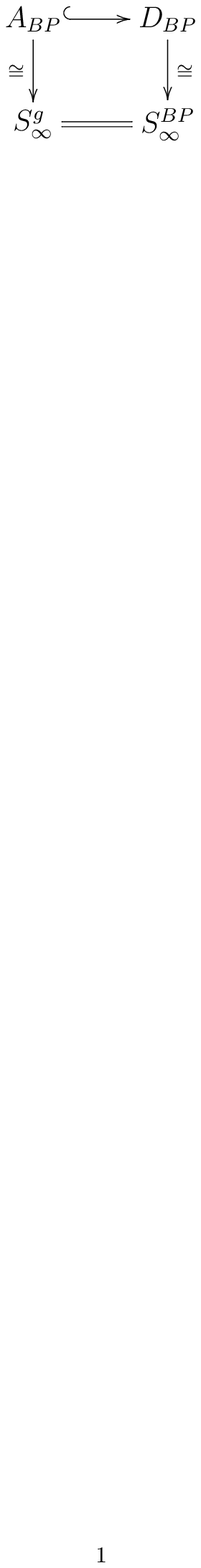}\\
\end{center}
\end{figure}
\vspace{-0.5cm}

\noindent and thus $D_{BP}=A_{BP}$.
\end{proof}

Now we need to explain how to prove the corresponding result for $BP$ when $p=2$ and
for $\MUp$ for all primes. The ideas involved are the same as those for $BP$ at
odd primes. Therefore we explain how the arguments need to be modified,
without giving all the details.

We begin by noting that, as in~\cite{Johnson1984} and~\cite{ccw}, we can characterize
$\kuz$ and $\kupz$ by systems of congruences. Thus, letting an operation in $\kuz$
act on
$\pi_{2i}(ku)$ as multiplication by $\lambda_i\in\Z$, we have
    $$
    \kuz\cong S_\infty^{ku}:=
    \left\{
    \lambda \in \prod_{k\geq 0} \Z \,\Big|\,
    \sum_{j=0}^i \frac{v_{i,j}}{d(i)}\lambda_j\in \Z\ \mbox{for all $i\geq 0$}
    \right\}.
    $$
Here $d(i)=\prod_{p\ \mbox{\scriptsize{prime}}} p^{\gamma_p(i)}$, where 
$\gamma_p(i)=\delta_p\left(\floor{i/(p-1)} \right)$. Descriptions of the integers $v_{i,j}$  can be found
in~\cite{Johnson1984} and~\cite{ccw}. (More precisely, these references give algorithms
for the $v_{i,j}$, but not a global formula.) We also have
    $$
    \kupz\cong S_\infty^{\kup}:=
    \left\{
    \lambda \in \prod_{k\geq 0} \Zpl \,\Big|\,
    \sum_{j=0}^i \frac{v_{i,j}}{d(i)}\lambda_j\in \Zpl\ \mbox{for all $i\geq 0$}
    \right\}.
    $$   
Using duality and~\cite[Props.\,3, 20]{ccw}, it is possible to write an equivalent
system of $p$-local $ku$ congruences explicitly, in terms of Gaussian polynomials, just as 
we did in the previous section for the Adams summand $g$. However, for our purposes
here, it is enough to note that we have one congruence for each non-negative integer $n$ and this takes the
form $p^{-\gamma_p(n)}(\lambda_n+y)\in \Zpl$, where $y$ is a $\Zpl$ 
linear combination
of $\lambda_0, \lambda_1, \dots, \lambda_{n-1}$. Of course, we may also define
$S_n^{ku}$ and $S_n^{\kup}$ in the obvious way, by imposing the relevant
congruences for $0\leq i\leq n$.

\begin{theorem}\label{BPtwocentre}
Let $p=2$. The centre of the ring of 
degree zero $BP$ operations is the 
Adams subalgebra:
$$
Z\left(BP^0(BP)\right)=
A_{BP}.
$$
\end{theorem}

\begin{proof}
As before, the centre is the ring of diagonal operations, so we
need to show that $A_{BP}=D_{BP}$.
Using Proposition~\ref{kubpinject}, 
$A_{BP}\cong \kutz\cong S_{\infty}^{ku_{(2)}}$. 
Again $D_{BP}\cong S_\infty^{BP}$. For each $n\geq 0$,
the inclusion $S_n^{ku_{(2)}}\subseteq S_n^{BP}$ follows from
Proposition~\ref{kubpinject}. To establish the reverse inclusions,
it is necessary to produce a suitable
$BP$ congruence for each non-negative integer $n$, of the form described
above for the congruences defining $S_{\infty}^{ku_{(2)}}$. 
Noting that, for $p=2$, $\gamma_p(n)=\delta_p(n)$, we see that
such a congruence is provided by
Proposition~\ref{specialcong} (which does hold for $p=2$). 
Just as in the previous section, we conclude that 
$S_n^{BP}=S_n^{ku_{(2)}}$ for all $n\geq 0$, and therefore
$S_{\infty}^{BP}=
S_{\infty}^{ku_{(2)}}$. The result follows.
\end{proof}

Now we turn to the results for $\MUp$.
We have $MU_*=\Z[x_1, x_2, \dots]$, where $|x_i|=2i$, and 
$MU_*(MU)=MU_*[b_1, b_2, \dots]$, where $|b_i|=2i$. Now let
$\theta\in D_{MU}$ act on $\pi_{2i}(MU)$ as multiplication by $\lambda_i\in\Z$.
Just as for $BP$, we can identify $D_{MU}$ with a set of sequences characterized
by a family of congruences. Thus we set $\thetab(b^\gamma)=\sum_{\delta}\Delta_\delta^\gamma(\lambda) x^\delta$,
and the $\Delta_\delta^\gamma(\lambda)$ are determined recursively
as finite rational linear combinations of the $\lambda_i$ 
by $\thetab\eta_R=\theta_*=\lambda$. Then we define      
    $$
    S_\infty^{MU}=\left\{\lambda\in\prod_{k\geq 0} \Z\,\Big|\, 
    \Delta_\delta^\gamma(\lambda)\in\Z 
    \mbox{ for all $\gamma, \delta$} \right\},
    $$
and we have $D_{MU}\cong S_\infty^{MU}$.
Similarly, we let
    $$
    S_\infty^{\MUp}=\left\{\lambda\in \prod_{k\geq 0} \Zpl\,\Big|\, \Delta_\delta^\gamma(\lambda)\in\Zpl 
    \mbox{ for all $\gamma, \delta$} \right\},
    $$
and we see that $D_{\MUp}\cong S_\infty^{\MUp}$.

\begin{theorem}\label{MUpcentre}
For any prime $p$, the centre of the ring of degree zero 
$\MUp$ operations is the Adams subalgebra:
$$
Z\left(\MUpz\right)=
A_{\MUp}.
$$
\end{theorem}

\begin{proof}
By Proposition~\ref{diag}, 
the centre is the ring of diagonal operations and we
need to show that $A_{\MUp}=D_{\MUp}$.
Using Propositions~\ref{kinject} and~\ref{ktinject}, for all primes $p$,
we have 
    $$
    A_{\MUp}\cong \kupz \cong S_\infty^{\kup}.
    $$
Also, $D_{\MUp}\cong S_\infty^{\MUp}$. Again the isomorphisms are given by sending 
operations to their actions on coefficients. So we need 
to show that $S_\infty^{\MUp}=S_\infty^{\kup}$. We will explain how
to exploit the results for $BP$ to show this.

Suppose that $\theta\in D_{\MUp}$ acts on $\pi_{2n}(\MUp)$ 
as multiplication by $\lambda_n$.
Using the same arguments as for $BP$,
it will be enough to find a particular family of $\MUp$ congruences
satisfied by the $\lambda_i$. 
Specifically, we need,
for each $n\geq 0$, an $\MUp$ congruence of the
form $p^{-\gamma_p(n)}(\lambda_n+y)\in \Zpl$, where $y$ is a $\Zpl$ 
linear combination
of $\lambda_0, \lambda_1, \dots, \lambda_{n-1}$. 
This is because $S_\infty^{\kup}$ is given by such a system of
congruences, as explained above.

Now ${\MUp}_*(\MUp)$ is a polynomial extension of $BP_*(BP)$ and we
write $i:BP_*(BP)\hookrightarrow {\MUp}_*(\MUp)$ for the inclusion. 
The elements $d_k$ of $BP_*(BP)$, constructed
in Proposition~\ref{specialcong}, give us elements $i(d_k)$ in ${\MUp}_*(\MUp)$.
Let $V_\lambda:{\MUp}_*(\MUp)\to \Zpl$ be given by
composing $\thetab:{\MUp}_*(\MUp)\to {\MUp}_*$ 
with the algebra map ${\MUp}_*\to \Zpl$ which sends $x_{p-1}$ to $1$
and $x_j$ to zero for all $j\neq p-1$.
Then, for $n=k(p-1)$ we 
can deduce the required congruence easily from the corresponding
$BP$ one, using $V_{\lambda}(i(d_k))=V_\mu(d_k)\in \Zpl$, 
where $\mu=(\lambda_0, \lambda_{p-1}, \lambda_{2(p-1)}, \dots)$.

For $p=2$, this gives all the required congruences. Finally,
let $p$ be an odd prime and let $n=k(p-1)+j$, with $0\leq j\leq p-2$.
Since $\gamma_p(n)=\delta_p(k)$, we need a congruence
of the form $p^{-\delta_p(k)}\left(\lambda_{n}+y\right)\in\Zpl$,
where $y$ is some $\Zpl$ linear combination of 
$\lambda_0,\dots, \lambda_{n-1}$.
By definition of $V_\lambda$, the sum over $r$
of the coefficients of $x_{p-1}^r$ in $\thetab(i(d_k))$
is $V_\lambda(i(d_k))$.
Now we have $\thetab(2 b_1)=(\lambda_1-\lambda_0)x_1$
and it follows that 
$\thetab(2^j b_1^j)=\left(\sum_{l=0}^j (-1)^{j-l}\binom{j}{l}\lambda_l\right)x_1^j$.
Using this, it may be checked that a congruence of the required type
is provided by the sum over $r$ of the coefficients of the 
terms $x_1^jx_{p-1}^r$ in $\thetab\left(2^jb_1^ji(d_k)\right)$. 
\end{proof}

We end with a result about integral $MU$. 

\begin{theorem}
$Z\left(MU^0(MU)\right)\cong \kuz$.
\end{theorem}

\begin{proof}
By Proposition~\ref{diag},
$Z\left(MU^0(MU)\right)=D_{MU}$. 
Also $D_{MU}\cong S_\infty^{MU}$ and $\kuz\cong S_\infty^{ku}$.
Now we have 
    $$
    S_\infty^{MU}=\bigcap_{p} S_\infty^{MU_{(p)}}\quad\mbox{and}\quad 
    S_\infty^{ku}=\bigcap_p S_\infty^{{ku}_{(p)}}.
    $$ 
The proof of 
Theorem~\ref{MUpcentre} shows that, for every prime $p$, 
$S_\infty^{{MU}_{(p)}}=S_\infty^{ku_{(p)}}$
and so $S_\infty^{MU}=S_\infty^{ku}$.
\end{proof}

\end{document}